\documentclass[preprint,10pt]{elsarticle}

\usepackage{amsmath,amssymb,amsbsy,amsfonts}
\usepackage{graphicx}
\usepackage{float}
\usepackage{morefloats}
\usepackage{multirow}
\usepackage{nomencl}
\usepackage{lipsum}
\usepackage{enumerate}
\newtheorem{thm}{Theorem}
\newtheorem{prop}[thm]{Proposition}

\newdefinition{rmk}{Remark}
\newproof{pf}{Proof}

\journal{Journal}
\begin{document}
\begin{frontmatter}
\title{An entropy admissible time splitting scheme for a conservation law model of manufacturing system}
\author{Tanmay Sarkar\corref{cor1}}
\address{Department  of  Mathematics,
Indian Institute of Technology Madras, \\Chennai-600 036, India}
\ead{tanmaysemilo@gmail.com}
\cortext[cor1]{Corresponding author. Tel: +91 9962024570}
\begin{abstract}
This paper deals with a splitting method applied to a conservation law model of manufacturing system incorporating yield loss. A splitting scheme has been proposed. The yield loss term is treated by solving implicitly an ordinary differential equation and the hyperbolic part is approximated by a finite volume scheme. Bounded variation stability has been studied. Due to yield loss, proposed scheme is total variation bounded. The convergence of the numerical solution towards entropy solution (in the Kruzkov sense) is proved. Numerical experiments are presented to demonstrate the performance of the scheme.
\end{abstract}
\begin{keyword}
hyperbolic conservation laws, splitting method, yield loss, BV stability, manufacturing system.
\end{keyword}
\end{frontmatter}
\numberwithin{equation}{section}

\section{Introduction}
Most problems of scientific interest are nonhomogeneous in nature. Dynamics of these problems are represented by hyperbolic conservation laws with source terms. In general, the appearance of the source term is either due to physical effects like exterior forces, release of mass or energy, chemical reacting gas etc. or due to geometrical effects like axisymmetric or cylindrical problems.

The purpose of this paper is to study an operator splitting procedure applied to a hyperbolic conservation law model of a manufacturing system. The model of manufacturing system has been introduced by Armbruster et al. in \cite{armbruster2006continuum}. Thereafter, it has been studied by several authors \cite{gottlich2005network, kirchner2006optimal, yield2014conservation, shang2011analysis}. Incorporating physical effects like yield loss, the conservation law model becomes nonhomogeneous. The model of manufacturing system is studied here as follows:
\begin{equation}\label{model_manuf}
\partial_t u(x,t)+\partial_x f(u)+y_l(x,t,u)=0,~~~t>0,~0< x\leq 1,
\end{equation}
where $u(x,t)$ is the density of the material at stage $x$ and time $t$. The flux function $f$ is given by
\begin{align*}
f(u)=v(WIP(t))u(x,t),~~~WIP(t)=\int_0^1u(x,t)dx.
\end{align*}
The term $y_l(x,t,u)$ represents the yield loss during the process. We assume that the velocity function $v$ is continuously differentiable. In the manufacturing system, with a given initial data
\begin{equation*}
u(x,0)=u_0(x),~~~~~0\leq x\leq 1,
\end{equation*}
natural input is the influx given as
\begin{equation*}
u(0,t)v(WIP(t))=\lambda(t),~~~t\geq 0.
\end{equation*}
Motivated by applications, we observe that the main objective of the model of manufacturing system is to analyze density distribution and outflux of the system. The outflux is given by $w(t)=v(WIP(t))u(1,t)$.
In this context, the homogeneous problem has been studied theoretically and numerically by several authors \cite{cutolo2011upwind, la2010control, sarkar2013_conservation}. The extensive theory of nonhomogeneous scalar can be found in \cite{kruvzkov1970}.

There have been significant contributions for the approximation of conservation laws involving source term in Chalabi \cite{chalabi1992_stable} and Schroll \cite{schroll1996finite}. These numerical schemes are based on explicit schemes. The solution of nonhomogeneous equation does not possess total variation diminishing (TVD) property because of the effect of source term which assists to increase the total variation. Sweby \cite{sweby89tvd} proposed a method based on the transformation of the dependent variable to reduce the nonhomogeneous scalar conservation law to a homogeneous one which possess the TVD property. It is well known that the explicit schemes are not appropriate for the numerical treatment of the source terms in several cases, this motivates us to use time splitting schemes.

Taking into account of source term $g=g(u)$ as a well behaved smooth function, several authors have investigated in this direction, one can refer in Crandall et al. \cite{crandall1980method} and Month$\acute{e}$ \cite{monthe01study}. Along this direction, error bound has been given by Tang and Teng in \cite{tang1995error} and convergence analysis has been studied by Langseth, Tveito and Winther in \cite{langseth1996convergence}. The objective of this paper is to relax the conditions on source term while taking into account of the dependence of space and time.

In this paper, we study the stability and convergence of the approximated solution obtained by splitting scheme where the ordinary differential equation part is handled by implicit scheme and the hyperbolic part is approximated using a finite volume monotone scheme. Considering the implicit character in source term, the proposed scheme is total variation bounded (TVB) and at the limit, satisfies entropy condition in the Kruzkov sense.

The paper is structured as follows. We start with theoretical investigation. For theoretical study, we contemplate the model (\ref{model_manuf}) in more general set up. In section 2, we present some preliminaries related to the nonhomogeneous scalar conservation laws. Section 3 concerns the splitting scheme where the stability estimate and convergence of the numerical solution towards the entropy solution is proved. Section 5 is devoted to the numerical investigation and subsequent discussion.
\section{Preliminaries}
The nonhomogeneous scalar conservation law to be investigated, is represented in this section by the following Cauchy problem:
\begin{equation}\label{nonhom1}
\partial_tu(x,t)+\partial_xf(u)=g(x,t,u),
\end{equation}
for $(x,t)\in\mathbb{R}\times[0,T];~T>0$ and
\begin{equation}\label{intial}
u(x,0)=u_0(x),~x\in\mathbb{R}
\end{equation}
with $u_0\in BV(\mathbb{R})\cap L^{1}(\mathbb{R})$. $BV(\mathbb{R})$ denotes the subspace of $L^1_{loc}$ consisting of functions with bounded variation, i.e.,
\begin{equation*}
BV(\mathbb{R})=\{v\in L^1_{loc}:\text{T.V.}(v)<\infty\},
\end{equation*}
where
\begin{equation*}
\text{T.V.}(v)=\sup_{h\neq 0}\int_{\mathbb{R}}\left|\displaystyle\frac{v(x+h)-v(x)}{h}\right|dx.
\end{equation*}
We assume that the function $g(x,t,u)$ satisfies the following properties:
\begin{enumerate}[(i)]
\item $g(x,t,u)$ is bounded for each fixed $u$ and continuous in t,
\item $|g(x,t,u_1)-g(x,t,u_2)|\leq L|u_1-u_2|$, for all $u_1$ and $u_2$ and $L$ is a constant independent of $x$ and $t$,
\item $\text{T.V.}(g(.,t,u))\leq B(t)$, for a bounded function $B(t)$ in $L^1[0,T]$.
\end{enumerate}
A bounded measurable function, $u(x,t)$, is a weak solution of (\ref{nonhom1}) and (\ref{intial}) if for all
$\phi\in C^1(\mathbb{R}\times[0,T])$ with compact support in $\mathbb{R}\times[0,T]$,
\begin{equation}\begin{split}
\displaystyle\int_0^T\int_{\mathbb{R}}(u(x,t)\phi_t(x,t)+f(u)\phi_x)(x,t)dxdt
+\int_{\mathbb{R}}u_0(x)\phi(x,0)dx\\
-\int_{\mathbb{R}}u(x,T)\phi(x,T)dx+\int_0^T\int_{\mathbb{R}}g(x,t,u)\phi(x,t) dxdt=0.
\end{split}\end{equation}
Since weak solutions are not uniquely determined in general by their initial data and additional principles, one needs to add an entropy condition to select the physically correct solution. We define the entropy condition in Kruzkov sense.

A bounded measurable function $u=u(x,t)$ is called an entropy solution of (\ref{nonhom1}) and (\ref{intial}) in $\mathbb{R}\times[0,T]$ if for any constant $k\in\mathbb{R}$ and any smooth function $\phi(x,t)\geq0$ with compact support in $\mathbb{R}\times[0,T]$, the following holds:
\begin{equation}\label{entropy}
\begin{split}
\displaystyle\int_0^T\int_{\mathbb{R}}\big(|u(x,t)-k|\phi_t+\text{sign}(u(x,t)-k)(f(u)-f(k))\phi_x\big)dxdt\\
+\int_{\mathbb{R}}\big(|u_0(x)-k|\phi(x,0)-|u(x,T)-k|\phi(x,T)\big)dx\\
+\displaystyle\int_0^T\int_{\mathbb{R}}\text{sign}(u(x,t)-k)g(x,t,u)\phi(x,t) dxdt\geq 0.
\end{split}\end{equation}
We observe that along the characteristic curves the solution $u(x,t)$ is not necessarily constant. For theoretical aspects, one can refer Kruzkov's result in \cite{kruvzkov1970}.

The spatial domain is divided into cells $I_j=[x_{j}-\frac{\Delta x}{2},x_{j}+\frac{\Delta x}{2})$ with centers at the point $x_j=j\Delta x,~j\in\mathbb{Z}$. Similarly, the time domain $[0,T]$ is discretized by $t^n=n\Delta t$ for $n\in\mathbb{N}$. Time strip is denoted by $J^n=[t^n,t^{n+1})$. Let $\chi_j^n$ be the characteristic function for the rectangle $R_j^n=I_j\times J^n$.

A weak solution $u(x,t)$ of (\ref{nonhom1})-(\ref{intial}) is approximated by a function $u^{\Delta}(x,t)$
defined on $\mathbb{R}\times[0,T]$ by
\begin{equation}
u^{\Delta}(x,t)=\displaystyle\sum_{n\in\mathbb{N}}\sum_{j\in\mathbb{Z}}\chi_j^n(x,t)u_j^n.
\end{equation}
The initial data is projected onto the space of piecewise constant functions by
\begin{equation*}
u^{\Delta}(x,0)=\sum_j \chi_j(x)u_j^0,~~~u_j^0=\frac{1}{\Delta x}\int_{I_j}u_0(x)dx,
\end{equation*}
where $\chi_{j}(x)$ is the characteristic function of the space grid $I_{j}$ for $j\in\mathbb{Z}$.

\section{Study on splitting scheme}
To take into account of nonhomogeneous character in the numerical solution of (\ref{nonhom1})-(\ref{intial}), we construct a splitting scheme using piecewise stationary data. The source term is handled by solving implicitly an ordinary differential equation and then treating the homogeneous part explicitly. We use the following discretized scheme:
\begin{align}\label{split_scheme1}
\displaystyle\frac{\bar{u}_j^n-u_j^n}{\Delta t}&=g(x_j,t^n,\bar{u}_j^n),
\end{align}
\begin{align}\label{split_scheme2}
u_j^{n+1}&=\bar{u}_j^n-\frac{\Delta t}{\Delta x}
[F(\bar{u}_j^n,\bar{u}_j^{n+1})-F(\bar{u}_j^{n-1},\bar{u}_j^n)],
\end{align}
where $F$ is the numerical flux, satisfies the following assumptions:
\begin{enumerate}[(i)]
\item $F$ is locally Lipschitz continuous function from $\mathbb{R}^2$ to $\mathbb{R}$,
\item $F(s,s)=f(s)$, i.e., numerical flux is consistent with the original flux,
\item $(a,b)\mapsto F(a,b)$, is non-decreasing with respect to $a$ and non-increasing with respect to $b$.
\end{enumerate}
The monotone flux scheme (\ref{split_scheme2}) can be written in the form:
\begin{equation*}
u_j^{n+1}=H(\bar{u}_{j-1}^n,\bar{u}_j^n,\bar{u}_{j+1}^n).
\end{equation*}
With the Courant Friedrichs Lewy (CFL) condition
\begin{equation}\label{CFL}
\frac{\Delta t}{\Delta x}\max_{p,q}\left[|F(u,q)-F(v,q)|+|F(p,u)-F(p,v)|\right]\leq |u-v|,
\end{equation}
for all $u,v\in S$, where
$S=\{w\in L^{\infty}(\mathbb{R}):\|w\|_{L^{\infty}({\mathbb{R}})}\leq A\|u_0\|_{L^{\infty}({\mathbb{R}})}\}$, $A$ is a positive constant, one can observe that monotone flux scheme (\ref{split_scheme2}) is monotone since
\begin{align*}
\frac{\partial H}{\partial \bar{u}_{j-1}^n}&=
\frac{\Delta t}{\Delta x}F_a(a,b)|_{(\bar{u}_{j-1}^n,\bar{u}_{j}^n)}\geq 0,\\
\frac{\partial H}{\partial \bar{u}_{j}^n}&=1-\frac{\Delta t}{\Delta x}
\big[F_a(a,b)|_{(\bar{u}_{j}^n,\bar{u}_{j+1}^n)}-F_b(a,b)|_{(\bar{u}_{j-1}^n,\bar{u}_{j}^n)}\big],
\end{align*}
using CFL condition, one can obtain $\displaystyle\frac{\partial H}{\partial \bar{u}_{j}^n}\geq 0$ and
\begin{align*}
\frac{\partial H}{\partial \bar{u}_{j+1}^n}&=-\frac{\Delta t}{\Delta x}
F_b(a,b)|_{(\bar{u}_{j}^n,\bar{u}_{j+1}^n)}\geq 0.
\end{align*}
\subsection{Stability Estimates}
\begin{prop}
Let $g(x_j,t^n,\bar{u}_j^n)$ satisfies the properties $(i)-(iii)$ in section 2. Then there exists a constant $L_g$ such that
\begin{equation}
|g(x_j,t^n,\bar{u}_j^n)|\leq L_g(1+ |\bar{u}_j^n|).
\end{equation}
\end{prop}
\begin{pf}
From properties $(ii)$, we have
\begin{equation*}
|g(x_j,t^n,\bar{u}_j^n)-g(x_j,t^n,\bar{v}_j^n)|\leq L|\bar{u}_j^n-\bar{v}_j^n|.
\end{equation*}
Let us choose $\bar{v}_j^n=0$. Then we have
\begin{equation*}
|g(x_j,t^n,\bar{u}_j^n)|\leq |g(x_j,t^n,0)|+L|\bar{u}_j^n|.
\end{equation*}
Since $g(x,t,u)$ is bounded for each fixed $u$, we obtain
\begin{equation*}
|g(x_j,t^n,\bar{u}_j^n)|\leq L_g(1+ |\bar{u}_j^n|),
\end{equation*}
for some constant $L_g$.
\end{pf}
\begin{prop}\label{tv_bound}
Let $g(x_j,t^n,\bar{u}_j^n)$ satisfies the properties $(i)-(iii)$ in section 2.
If the CFL condition (\ref{CFL}) and $L_g\Delta t<1$ are satisfied, then the schemes (\ref{split_scheme1})-(\ref{split_scheme2}) satisfy the following estimates:
\begin{align}
\|u^{n+1}\|_{L^{\infty}(\mathbb{Z})}\leq \exp(C_0T)\|u_{0}\|_{L^{\infty}(\mathbb{Z})}
\end{align}
\begin{align}
\text{T.V.}(u^{n+1})\leq \exp(C_0T)(\text{T.V.}(u_{0})+\|B\|_1),
\end{align}
where $C_0$ is a positive constant.
\end{prop}
\begin{pf}
The scheme (\ref{split_scheme1}) can be written as
\begin{align*}
\bar{u}_j^n &= u_j^n-\Delta t g(x_j,t^n,\bar{u}_j^n),\\
|\bar{u}_j^n|& \leq |u_j^n|+\Delta t |g(x_j,t^n,\bar{u}_j^n)|\leq |u_j^n|+\Delta t L_g(1+|\bar{u}_j^n|),\\
(1-& \Delta t L_g)|\bar{u}_j^n| \leq |u_j^n|+\Delta t L_g,\\
|\bar{u}_j^n|& \leq \frac{1}{1-\Delta t L_g}|u_j^n|+\frac{\Delta t L_g}{1-\Delta t L_g}.
\end{align*}
If we set $C=\displaystyle\frac{L_g}{1-\Delta t L_g}$, the above inequalities can be written as
\begin{align*}
|\bar{u}_j^n|& \leq (1+ C\Delta t)|u_j^n|+C\Delta t.
\end{align*}
We assert the following:
\begin{align*}
|\bar{u}_j^n|& \leq (1+ C_1\Delta t)|u_j^n|,
\end{align*}
for some constant $C_1$ depending on $\Delta t$. This implies
\begin{equation}\label{condition1}
\|\bar{u}^n\|_{L^{\infty}(\mathbb{Z})} \leq (1+ C_1\Delta t)\|u^n\|_{L^{\infty}(\mathbb{Z})}.
\end{equation}
The scheme (\ref{split_scheme2}) can be written as the following form:
\begin{equation*}
u^{n+1}_j=(1-d_{j+\frac{1}{2}}-e_{j-\frac{1}{2}})\bar{u}_j^n+d_{j+\frac{1}{2}}\bar{u}_{j+1}^n
+e_{j-\frac{1}{2}}\bar{u}_{j-1}^n,
\end{equation*}
where
\[ d_{j+\frac{1}{2}} = \left\{ \begin{array}{ll}
 \frac{k}{h}\frac{F(\bar{u}_j^n,\bar{u}_{j+1}^{n})-f(\bar{u}_j^n)}{\bar{u}_j^n-\bar{u}_{j+1}^{n}},
 & \mbox{if  $\bar{u}_j^n \neq \bar{u}_{j+1}^{n})$}\\
  0,                & \mbox{if  $\bar{u}_j^n = \bar{u}_{j+1}^{n})$}.\end{array} \right. \]
 and
 \[ e_{j-\frac{1}{2}} = \left\{ \begin{array}{ll}
 \frac{k}{h}\frac{F(\bar{u}_{j-1}^n,\bar{u}_{j}^{n})-f(\bar{u}_j^n)}{\bar{u}_{j-1}^n-\bar{u}_{j}^{n}},
 & \mbox{if  $\bar{u}_j^n \neq \bar{u}_{j-1}^{n})$}\\
  0,                & \mbox{if  $\bar{u}_j^n = \bar{u}_{j-1}^{n})$}.\end{array} \right. \]
  Thanks to the monotonicity of $F$ and $f(\bar{u}_j^n)=F(\bar{u}_j^n,\bar{u}_j^n)$, we can conclude that
  \begin{align}\label{condition2}
  \|{u}^{n+1}\|_{L^{\infty}(\mathbb{Z})} \leq \|\bar{u}^n\|_{L^{\infty}(\mathbb{Z})}
  \end{align}
 Conditions (\ref{condition1}) and (\ref{condition2}) assert that
 \begin{equation*}
 \|{u}^{n+1}\|_{L^{\infty}(\mathbb{Z})} \leq (1+ C_1\Delta t)\|u^n\|_{L^{\infty}(\mathbb{Z})}
 \end{equation*}
 which implies
 \begin{equation*}
 \|{u}^{n}\|_{L^{\infty}(\mathbb{Z})} \leq e^{C_1n\Delta t}\|{u}_0\|_{L^{\infty}(\mathbb{Z})}
 \leq e^{C_1T}\|{u}_{0}\|_{L^{\infty}(\mathbb{Z})}.
 \end{equation*}
 Using the similar arguments we get
 \begin{align*}
 |\bar{u}_{j+1}^n-\bar{u}_{j}^n|& \leq |u_{j+1}^n-u_{j}^n|
 +\Delta t|g(x_{j+1},t^n,\bar{u}_{j+1}^n)-g(x_{j},t^n,\bar{u}_{j}^n)|\\
 & \leq |u_{j+1}^n-u_{j}^n|
 +\Delta t|g(x_{j+1},t^n,\bar{u}_{j+1}^n)-g(x_{j},t^n,\bar{u}_{j+1}^n)|\\
 &+\Delta t|g(x_{j},t^n,\bar{u}_{j+1}^n)-g(x_{j},t^n,\bar{u}_{j}^n)|\\
 & \leq |u_{j+1}^n-u_{j}^n|+ \Delta t |B(t_n)|+ L\Delta t|\bar{u}_{j+1}^n-\bar{u}_{j}^n|\\
 (1-L\Delta t)& |\bar{u}_{j+1}^n-\bar{u}_{j}^n|\leq |u_{j+1}^n-u_{j}^n|+ \Delta t |B(t_n)|.
 \end{align*}
 The above inequality implies that
 \begin{align}\label{tv1}
 |\bar{u}_{j+1}^n-\bar{u}_{j}^n|& \leq (1+C\Delta t)[|u_{j+1}^n-u_{j}^n|+ \Delta t |B(t_n)|].
 \end{align}
 Taking into account of the CFL condition (\ref{CFL}) and given condition, we can easily show that
 \begin{equation*}
 \text{T.V.}(u^{n+1})\leq \text{T.V.}(\bar{u}^n).
 \end{equation*}
 Using the above inequality in (\ref{tv1}), we obtain
 \begin{align*}
 |{u}_{j+1}^{n+1}-{u}_{j}^{n+1}|& \leq (1+C\Delta t)[|u_{j+1}^n-u_{j}^n|+ \Delta t |B(t_n)|].
 \end{align*}
 Thus we have
 \begin{align*}
 \text{T.V.}(u^n)& \leq e^{Cn\Delta t}[\text{T.V.}(u_0)+\|B\|_1]
 \leq e^{CT}[\text{T.V.}(u_0)+\|B\|_1].
 \end{align*}
 Moreover, there exists $\Delta t_0>0$ such that
 \begin{align*}
 \forall~ \Delta t<\Delta t_0,~~~ C<\frac{L_g}{1-L_g\Delta t_0}=:C_0.
 \end{align*}
 Thus we obtain
 \begin{align*}
 \|{u}^{n}\|_{L^{\infty}(\mathbb{Z})} &\leq e^{C_0T}\|{u}_0\|_{L^{\infty}(\mathbb{Z})}\\
 \text{T.V.}(u^{n})&\leq e^{C_0T}\text{T.V.}(u_0)
 \end{align*}
\end{pf}
\begin{rmk}
We can also establish $BV$ Stability in time. Using the similar arguments as mentioned above and thanks to the property $g(x,t,u)$ is continuous in $t$, one can obtain the following:
\begin{equation*}
\|u_j\|_{L^{\infty}(\mathbb{N})}\leq M,\text{  and    } \text{T.V.}(u_j)\leq M,
\text{   for } j\in\mathbb{Z},
\end{equation*}
where $M$ is some constant.
\end{rmk}
\begin{prop}\label{discrete_entropy}
{\bf{(Discrete Entropy Inequality)}}
Let $u^{\Delta t}$ be the approximate solution of (\ref{nonhom1})-(\ref{intial}) using the splitting scheme (\ref{split_scheme1})-(\ref{split_scheme2}). Assume that the scheme is a monotone flux scheme and $g(x,t,u)$ satisfies the properties $(i)-(iii)$ in section 2. Under the CFL condition (\ref{CFL}), the following inequality holds:
\begin{equation}\begin{split}
(|u_j^{n+1}-k|-|u_j^n-k|)+\frac{\Delta t}{\Delta x}[F(\bar{u}_j^n\vee k,\bar{u}_{j+1}^n\vee k)-
F(\bar{u}_j^n\wedge k,\bar{u}_{j+1}^n\wedge k)\\
-F(\bar{u}_{j-1}^n\vee k,\bar{u}_{j}^n\vee k)+F(\bar{u}_{j-1}^n\wedge k,\bar{u}_{j}^n\wedge k)]
\leq \text{sign}(\bar{u}_{j}^n-k)\Delta t g(x_j,t^n,\bar{u}_{j}^n),\\~~~\forall~j\in\mathbb{Z},~\forall~n\in\mathbb{N}\text{ and }\forall~k\in\mathbb{R},
\end{split}\end{equation}
where $r_1\vee r_2~(resp.~r_1\wedge r_2)$ denotes the maximum (resp. minimum) of the two real numbers $r_1$ and $r_2$.
\end{prop}
\begin{pf}
Taking into account of CFL condition (\ref{CFL}) in scheme (\ref{split_scheme2}) and using the monotonicity properties of $F$,
\begin{equation*}
u_j^{n+1}=H(\bar{u}_{j-1}^n,\bar{u}_{j}^n,\bar{u}_{j+1}^n),~~~\forall~j\in\mathbb{Z},~\forall~n\in\mathbb{N},
\end{equation*}
where $H$ is a function from $\mathbb{R}^3$ to $\mathbb{R}$, we have shown
\begin{equation*}
\frac{\partial H}{\partial \bar{u}_{l}^n}\geq 0,~~~\text{for }l=j-1,j,j+1.
\end{equation*}
Also, we observe that $H(k,k,k)=k$, for all $k\in\mathbb{R}$. Thus, we have the following:
\begin{align*}
u_j^{n+1}&\leq H(\bar{u}_{j-1}^n\vee k,\bar{u}_{j}^n\vee k,\bar{u}_{j+1}^n\vee k),\\
k&\leq H(\bar{u}_{j-1}^n\vee k,\bar{u}_{j}^n\vee k,\bar{u}_{j+1}^n\vee k),
\end{align*}
for all $k\in\mathbb{R}$, which implies
\begin{align*}
u_j^{n+1}\vee k &\leq H(\bar{u}_{j-1}^n\vee k,\bar{u}_{j}^n\vee k,\bar{u}_{j+1}^n\vee k).
\end{align*}
Using similar arguments as above, we obtain the following
\begin{align*}
u_j^{n+1}\wedge k &\geq H(\bar{u}_{j-1}^n\wedge k,\bar{u}_{j}^n\wedge k,\bar{u}_{j+1}^n\wedge k).
\end{align*}
Using the above inequalities, we have
\begin{align*}
|u_j^{n+1}-k|&=(u_j^{n+1}\vee k)-(u_j^{n+1}\wedge k)\\
 & \leq H(\bar{u}_{j-1}^n\vee k,\bar{u}_{j}^n\vee k,\bar{u}_{j+1}^n\vee k)-H(\bar{u}_{j-1}^n\wedge k,\bar{u}_{j}^n\wedge k,\bar{u}_{j+1}^n\wedge k)\\
 &\leq (\bar{u}_{j}^n\vee k)-\frac{\Delta t}{\Delta x}[F(\bar{u}_j^n\vee k,\bar{u}_{j+1}^n\vee k)-F(\bar{u}_{j-1}^n\vee k,\bar{u}_{j}^n\vee k)]\\
 &-(\bar{u}_{j}^n\wedge k)+\frac{\Delta t}{\Delta x}[F(\bar{u}_{j}^n\wedge k,\bar{u}_{j+1}^n\wedge k)-F(\bar{u}_{j-1}^n\wedge k,\bar{u}_{j}^n\wedge k)]\\
 &\leq |\bar{u}_{j}^n-k|-\frac{\Delta t}{\Delta x}[F(\bar{u}_j^n\vee k,\bar{u}_{j+1}^n\vee k)-F(\bar{u}_{j-1}^n\vee k,\bar{u}_{j}^n\vee k)]\\
 & +\frac{\Delta t}{\Delta x}[F(\bar{u}_{j}^n\wedge k,\bar{u}_{j+1}^n\wedge k)-F(\bar{u}_{j-1}^n\wedge k,\bar{u}_{j}^n\wedge k)]
\end{align*}
Scheme (\ref{split_scheme1}) yields the following inequalities
\begin{align*}
|u_j^{n+1}-k|& \leq |{u}_{j}^n-k|+\text{sign}(\bar{u}_{j}^n-k)\Delta t g(x_j,t^n,\bar{u}_{j}^n)
-\frac{\Delta t}{\Delta x}[F(\bar{u}_j^n\vee k,\bar{u}_{j+1}^n\vee k)\\
& -F(\bar{u}_{j-1}^n\wedge k,\bar{u}_{j}^n\wedge k)-F(\bar{u}_{j}^n\wedge k,\bar{u}_{j+1}^n\wedge k)
+F(\bar{u}_{j-1}^n\wedge k,\bar{u}_{j}^n\wedge k)].
\end{align*}
That completes the proof of the proposition.
\end{pf}
\subsection{Convergence}
\begin{thm}
Let $g(x,t,u)$ satisfies the properties $(i)-(iii)$ in section 2. Let $u_0\in BV(\mathbb{R})\cap L^1(\mathbb{R})$ and $f$ is Lipschitz continuous. If the CFL condition (\ref{CFL}) and $L_g\Delta t<1$
are satisfied, then the approximate solution $u^{\Delta t}(x,t)$ constructed by the splitting scheme (\ref{split_scheme1})-(\ref{split_scheme2}) converges in $L^1_{loc}(\mathbb{R}\times[0,T])$ towards the unique entropy solution of (\ref{nonhom1})-(\ref{intial}) as $\Delta t$ tends to zero.
\end{thm}
\begin{pf}
\begin{enumerate}
\item From proposition \ref{tv_bound}, the sequence $u^{\Delta t}$ is bounded in $L^{\infty}(\mathbb{R}\times[0,T])\cap BV(\mathbb{R}\times[0,T])$. By Helley's theorem, there exists a subsequence still denoted by $u^{\Delta t}(x,t)$ converges to $u(x,t)$ in $L^1_{loc}(\mathbb{R}\times[0,T])$.\\
    Inequalities (\ref{condition1}) and (\ref{tv1}) yield the following:
    \begin{align*}
    \|\bar{u}^n\|_{L^{\infty}(\mathbb{Z})}& \leq \bar{C}\|u_0\|_{L^{\infty}(\mathbb{Z})},\\
    \text{T.V.}(\bar{u}^n)& \leq \bar{C}[\text{T.V.}(u_0)+\|B\|_1],
    \end{align*}
    for some constant $\bar{C}$. Now we define
    \begin{align*}
    \bar{u}^{\Delta t}(x,t)=\bar{u}_j^n,~~~for~(x,t)\in I_j\times J^n.
    \end{align*}
    Again by Helley's theorem, there exists a subsequence $\bar{u}^{\Delta t}(x,t)$ which converges towards $v(x,t)$ in $L^1_{loc}(\mathbb{R}\times[0,T])$.
\item Now we will show that $u(x,t)=v(x,t)$ in $L^1_{loc}(\mathbb{R}\times[0,T])$. Let $M$ be a compact set on $\mathbb{R}\times[0,T]$, then we have
    \begin{equation*}
    \|u-v\|_{L^1(M)}\leq \|u-u^{\Delta t}\|_{L^1(M)}+\|u^{\Delta t}-\bar{u}^{\Delta t}\|_{L^1(M)}
    +\|\bar{u}^{\Delta t}-v\|_{L^1(M)}.
    \end{equation*}
    We know that
    \begin{align*}
    |\bar{u}_j^n-u_j^n|&=|\Delta t g(x_j,t^n,\bar{u}_j^n)|
    & \leq \Delta t L_g(1+\bar{u}_j^n)
    & \leq \Delta t L_g(1+\bar{C}\|u_0\|_{L^{\infty}(\mathbb{Z})}),
    \end{align*}
    which implies
    \begin{align*}
    \|u^{\Delta t}-\bar{u}^{\Delta t}\|_{L^1(M)}\leq \text{mes}(M)\Delta t L_g(1+\bar{C}\|u_0\|_{L^{\infty}(\mathbb{Z})}).
    \end{align*}
    Hence $\|u^{\Delta t}-\bar{u}^{\Delta t}\|_{L^1(M)}\rightarrow 0$ as $\Delta t\rightarrow 0$.\\
    This proves that $u=v$ in $L^1_{loc}(\mathbb{R}\times[0,T])$. It remains to show that $u(x,t)$ satisfies entropy condition in the sense (\ref{entropy}).
\item We set $\eta(u)=|u-k|$ and $G(\bar{u},\bar{v})=F(\bar{u}\vee k,\bar{v}\vee k)-F(\bar{u}\wedge k,\bar{v}\wedge k)$.
    Since scheme (\ref{split_scheme2}) is monotone under the CFL condition (\ref{CFL}), by proposition \ref{discrete_entropy}, we have
    \begin{equation}\label{eta_entropy}
    \begin{split}
    (\eta(u_j^{n+1})-\eta(u_j^n))+\frac{\Delta t}{\Delta x}
    [G(\bar{u}_j^n,\bar{u}_{j+1}^n)-G(\bar{u}_{j-1}^n,\bar{u}_{j}^n)]\\
    \leq \eta'(\bar{u}_j^n)\Delta t g(x_j,t^n,\bar{u}_j^n).
    \end{split}\end{equation}
    Let $\phi$ be a nonnegative differentiable function on $\mathbb{R}\times[0,T]$. Multiply (\ref{eta_entropy}) by $\Delta x\phi (x_j,t^n)$ and taking summation over all $j\in\mathbb{Z}$ and $n\geq 0$, we obtain
    \begin{equation}\label{prf_entropy}
    S_1+S_2\leq S_3.
    \end{equation}
    Using discrete integration by parts
    \begin{align*}
    S_1&= \Delta x \sum_{n=0}^N\sum_{j\in\mathbb{Z}}(\eta(u_j^{n+1})-\eta(u_j^n))\phi_j^n\\
    &=\Delta x \Delta t \sum_{n=0}^N\sum_{j\in\mathbb{Z}}\frac{\eta(u_j^{n+1})-\eta(u_j^n)}
    {\Delta t}\phi_j^n\\
    &=\Delta x\left(-\sum_{j\in\mathbb{Z}}\eta(u_j^0)\phi_j^0-\Delta t\sum_{n=1}^N\sum_{j\in\mathbb{Z}}
    {\eta(u_j^n)}\frac{\phi_j^n-\phi_j^{n-1}}{\Delta t}
    +\sum_{j\in\mathbb{Z}}\eta(u_j^{N+1})\phi_j^N\right).
    \end{align*}
    Similarly, for $S_2$ and $S_3$ we get
    \begin{align*}
    S_2&=\Delta t\Delta x \sum_{n=0}^N\sum_{j\in\mathbb{Z}}\frac{1}{\Delta x}
    \left(G(\bar{u}_j^n,\bar{u}_{j+1}^n)-G(\bar{u}_{j-1}^n,\bar{u}_{j}^n)\right)\phi_j^n\\
    &=-\Delta t\Delta x \sum_{n=0}^N\sum_{j\in\mathbb{Z}}
    G(\bar{u}_j^n,\bar{u}_{j+1}^n)\frac{\phi_{j+1}^n-\phi_{j}^n}{\Delta x}\\
    S_3&=\Delta x \sum_{n=0}^N\sum_{j\in\mathbb{Z}}
    \eta'(\bar{u}_j^n)\Delta t g(x_j,t^n,\bar{u}_j^n)\phi_j^n\\
    &=\Delta t\Delta x\sum_{n=0}^N\sum_{j\in\mathbb{Z}}\frac{\eta(\bar{u}_j^n)-\eta(k)}{\bar{u}_j^n-k}
    g(x_j,t^n,\bar{u}_j^n)\phi_j^n.
    \end{align*}
    Now let $\Delta t\rightarrow 0$. It is reasonably straightforward, using the 1-norm convergence of $u^{\Delta t}$ to $u$ and smoothness of $\phi$, to show that as $\Delta t\rightarrow 0$
    \begin{align*}\begin{split}
    S_1\rightarrow -\int_{\mathbb{R}}\eta(u_0)\phi(x,0)dx-\int_0^T\int_{\mathbb{R}}\eta(u)\phi_t(x,t) dxdt\\
    +\int_{\mathbb{R}}\eta(u(x,T))\phi(x,T)dx.
    \end{split}\end{align*}
    The term $S_2$ needs additional description. We know
    \begin{equation*}
    G(\bar{u}_j^n,\bar{u}_{j+1}^n)=F(\bar{u}_j^n\vee k,\bar{u}_{j+1}^n\vee k)
    -F(\bar{u}_j^n\wedge k,\bar{u}_{j+1}^n\wedge k).
    \end{equation*}
    Since $F$ is consistent with $f$, we have
    \begin{align*}
    |F(\bar{u}_j^n,\bar{u}_{j+1}^n)-f(\bar{u}_j^n)|
    =|F(\bar{u}_j^n,\bar{u}_{j+1}^n)-F(\bar{u}_j^n,\bar{u}_j^n)|\leq \|F_b\|_{\infty}|\bar{u}_{j+1}^n-\bar{u}_{j}^n|.
    \end{align*}
    Thanks to the total bounded variation of $\bar{u}^{\Delta t}$, one can observe that numerical flux function can be approximated by $f(\bar{u}_j^n)$ with errors that vanish almost everywhere. Consequently, as $\Delta t\rightarrow 0$, we obtain the following:
    \begin{equation*}
    G(\bar{u}_j^n,\bar{u}_{j+1}^n)\rightarrow sign(u-k)(f(u)-f(k)).
    \end{equation*}
    This implies as $\Delta t\rightarrow 0$,
    \begin{equation*}
    S_2\rightarrow \int_0^T\int_{\mathbb{R}}\text{sign}(u-k)(f(u)-f(k))\phi_x(x,t) dxdt.
    \end{equation*}
    Using the properties $(i)-(iii)$ of $g(x,t,u)$ in section 2, we obtain
    \begin{equation*}
    \lim_{\Delta t\rightarrow 0}S_3=\int_0^T\int_{\mathbb{R}}\eta'(u) g(x,t,u)\phi(x,t) dxdt.
    \end{equation*}
    From (\ref{prf_entropy}), we can conclude that $u(x,t)$ satisfies the entropy condition in the sense (\ref{entropy}) as $\Delta t$ tends to zero.
\end{enumerate}
\end{pf}
\section{Numerical implementation and discussion}
In this section, we present some numerical experiments that demonstrate the performance of the proposed scheme. We focus on the model of manufacturing system and describe two cases regarding this. We use implicit scheme for the part involving source term and for the hyperbolic part, we consider finite volume method (for more details, refer \cite{leveque2002finite}).
\subsection{Test Case-I}
\begin{figure}[h!]
    \centering
    \includegraphics[width=12cm]{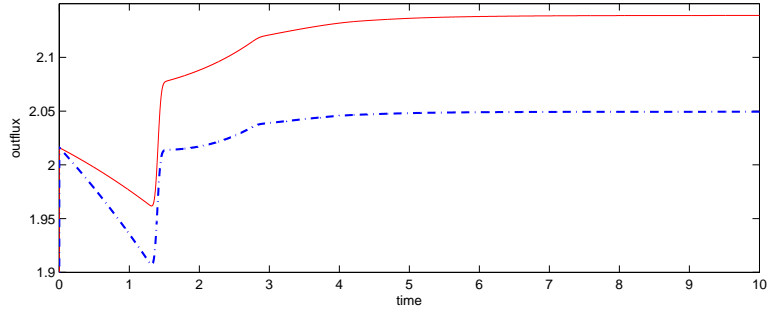}
    \caption{Armbruster et al. \cite{armbruster2006continuum} outflux (without yield loss) and outflux incorporating constant yield loss}
    \label{fig:a3constant}
\end{figure}
\begin{figure}[h!]
    \centering
    \includegraphics[width=9cm]{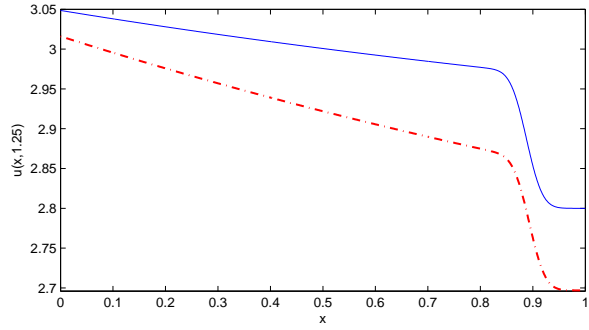}
    \caption{snapshots of the density $u(x)$ without yield loss (Armbruster et al. \cite{armbruster2006continuum}) and incorporating constant yield loss at time t=1.25}
    \label{fig:a3uT}
\end{figure}
\begin{figure}[h!]
    \centering
    \includegraphics[width=8cm]{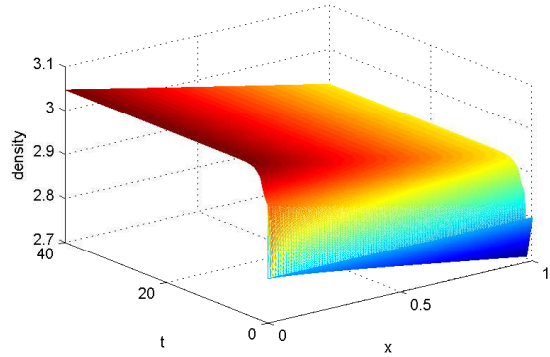}
    \caption{density distribution with constant yield loss}
    \label{fig:a3den}
\end{figure}
In the first case, we reproduce the result in Armbruster et al. \cite{armbruster2006continuum} as a part of validation of our scheme. The velocity term considered as follows:
\begin{align*}
v(WIP(t))&=v_0\left(1-\frac{WIP(t)}{L_m}\right),
\end{align*}
where $v_0$ is speed for the empty factory and $L_m$ represents the maximal load of the manufacturing system. Influx in the system is considered as
\[\lambda (t)= \left\{ \begin{array}{ll}
               2.016   & \mbox{$~~~~~$for  $t<0$}\\
               2.139   & \mbox{$~~~~~$for  $t>0$}.\end{array} \right. \]
At steady state, the density in the system is considered as 2.8. In the above set up, we carry out the simulation using the proposed scheme. We assume that there is $3\%$ yield loss in the density throughout the process.

Figure \ref{fig:a3constant} demonstrates the overview of the outflux. Outflux without yield loss and with yield loss have been exhibited. The outflux initially declines when the influx increases. This is an important observation in manufacturing system. As expected, the reduction in outflux is due to the yield loss. The outflux in both the cases become stable since there are no changes in influx and yield loss as time progresses.

Figure \ref{fig:a3uT} presents density distribution at $t=1.25$ in both cases. In the yield loss case, as we progress in space direction, reduction in density becomes higher due to the nonlinearity in the flux function. In both the cases, density is asymptotically approaching to a steady state. Figure \ref{fig:a3den} analyzes the overview of density distribution throughout the manufacturing system in the yield loss case.

\subsection{Test Case-II}
\begin{figure}[h!]
    \centering
    \includegraphics[width=7cm]{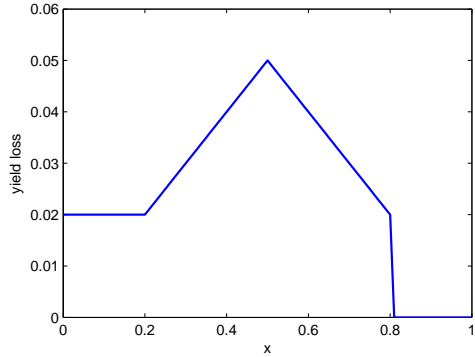}
    \caption{profile of yield loss during the process}
    \label{fig:yl_plot}
\end{figure}
\begin{figure}
    \centering
    \includegraphics[width=8cm]{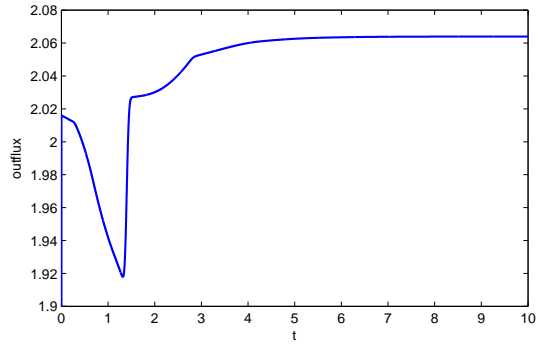}
    \caption{outflux in the manufacturing system during the process}
    \label{fig:outflux_yl}
\end{figure}
\begin{figure}[h!]
    \centering
    \includegraphics[width=8cm]{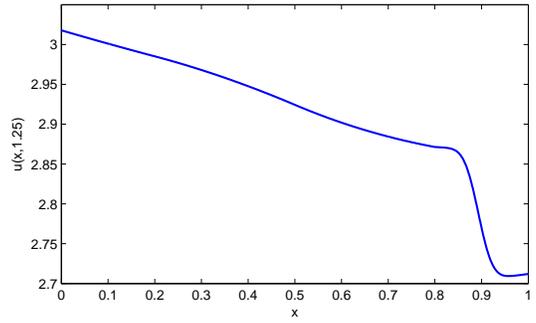}
    \caption{snapshot of the density at time $t=1.25$}
    \label{fig:yl_T}
\end{figure}
\begin{figure}[h!]
    \centering
    \includegraphics[width=8cm]{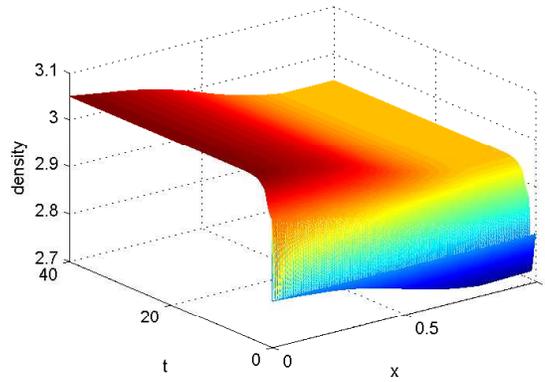}
    \caption{density distribution of the manufacturing system during the process}
    \label{yl_noncon}
\end{figure}
In this case, we consider yield loss as a piecewise linear function which depends on the space variable. Figure \ref{fig:yl_plot} illustrates the profile of yield loss during the process. All the other parameters remain same as in Test Case-I. Figure \ref{fig:outflux_yl} asserts that the outflux reacts a bit slower than constant yield loss case. Since yield loss does not explicitly depend in time and influx remains constant as time progresses, outflux becomes stable. In Figure \ref{fig:yl_T}, we observe that the asymptotic approach of the density to a stable state is not linear due to the yield loss profile. Figure \ref{yl_noncon} provides the density distribution in a manufacturing system. The effect of yield loss in density can be visualized without much difficulty.

\subsection*{Acknowledgement}
Author would like to thank Council of Scientific and Industrial Research (CSIR), New Delhi for the financial support (Ref. no. 09/084 (0505)/2009-EMR-I).

\bibliographystyle{model3a-num-names}
\bibliography{plain_split}
\end{document}